\documentclass[a4paper,12pt]{article}
\usepackage{longtable}
\title{Remarks on abelian surfaces\\ in nonsingular toric Fano 4-folds}

\date{}

\author{{\sc Hiroshi Sato}\thanks{Partly supported by the Grant-in-Aid for JSPS Fellows, The Ministry of Education, Science, Sports and Culture, Japan.
\newline
\hspace*{1.5em} {\em $2000$ Mathematics Subject Classification\/}.
Primary 14M25;
Secondary 14J45, 14K12.}}

\newtheorem{Thm}{Theorem}[section]
\newtheorem{Prop}[Thm]{Proposition}
\newtheorem{Cor}[Thm]{Corollary}
\newtheorem{Conj}[Thm]{Conjecture}
\newtheorem{Lem}[Thm]{Lemma}

\newtheorem{Prob}[Thm]{Problem}
\newtheorem{Def}[Thm]{Definition}
\newtheorem{Rem}[Thm]{Remark}
\newtheorem{Ex}[Thm]{Example}
\newtheorem{Not}[Thm]{Notation}

\newcommand{\proof}{Proof. \quad}
\newcommand{\qed}{\hfill q.e.d.}

\newcommand{\Hom}{\mathop{\rm Hom}\nolimits}

\newcommand{\G}{\mathop{\rm G}\nolimits}

\begin{document}
\maketitle

\begin{abstract}                                                   

In this paper, we investigate whether the 124 nonsingular toric Fano 4-folds admit totally nondegenerate embeddings from abelian surfaces or not. As a result, we determine the possibilities for such embeddings except for the remaining 21 nonsingular toric Fano 4-folds.

\end{abstract}

\section{Introduction}\label{intro}

\thispagestyle{empty}

\hspace{5mm} There exist no embeddings from abelian surfaces into nonsingular projective toric $3$-folds over ${\bf C}$ (see, e.g., Kajiwara \cite{kajiwara2} and \cite{kajiwara1}). So the next problem is to study which nonsingular projective toric $4$-folds admit embeddings from abelian surfaces. This problem was considered by many people (see Horrocks-Mumford \cite{horrocks1}, Hulek \cite{hulek1}, Kajiwara \cite{kajiwara2}, \cite{kajiwara1}, Lange \cite{lange1} and Sankaran \cite{sankaran1}). In this paper, we consider the following problem.

\begin{Prob}\label{mainfanoprob}
{\rm
Which nonsingular toric Fano $4$-fold admits a {\em totally nondegenerate embedding} from an abelian surface (see Definition \ref{defoftndabel})?
}

\end{Prob}

There exist exactly $124$ nonsingular toric Fano $4$-folds up to isomorphism (see Batyrev \cite{batyrev4} and Sato \cite{sato1}). We give a partial answer to Problem \ref{mainfanoprob} (see Theorem \ref{mt}).

The content of this paper is as follows: In Section \ref{prepare}, we recall the definition of a {\em totally nondegenerate embedding}. In Section \ref{somecri}, we describe criteria for the non-existence of totally nondegenerate finite morphisms, and using these criteria, we show the non-existence for some nonsingular projective toric $4$-folds. In Section \ref{main2blowup}, we consider the relationship between $2$-blow-ups of toric $4$-folds and totally nondegenerate finite morphisms. As a result, we can derive the main result in Section \ref{mainresult}. In Section \ref{secex}, we show the non-existence of totally nondegenerate finite morphisms for some nonsingular toric Fano $4$-folds. In Section \ref{mainresult}, we obtain the main result.

The author wishes to thank Professors Tadao Oda, Masa-Nori Ishida and Takeshi Kajiwara for advice and encouragement.

\section{totally nondegenerate embedding}\label{prepare}

\hspace{5mm} The following notation is used throughout this paper. For fundamental properties of the toric geometry, see Oda \cite{oda2}.

Let $N:={\bf Z}^{4}$ and $M:=\Hom_{\bf Z}(N,{\bf Z})$ the dual group. For a finite complete nonsingular fan $\Sigma$ in $N$ and $0\leq i\leq 4$, we put $\Sigma(i):=\left\{ \sigma\in\Sigma\; |\;\dim\sigma=i\right\}$. Each $\tau\in\Sigma(1)$ determines a unique element $e(\tau)\in N$ which generates the semigroup $\tau\cap N$. We call
$$\G(\Sigma):=\left\{e(\tau)\in N\; |\;\tau\in\Sigma(1)\right\}$$
the set of primitive generators of $\Sigma$. Let $X$ be the complete nonsingular toric $4$-fold corresponding to $\Sigma$. Let $\G(\Sigma)=\left\{ x_1,\ldots,x_n\right\}$ and let $\left\{ D_1,\ldots,D_n\right\}$ be the corresponding $T_N$-invariant prime divisors on $X$. Especially, the Picard number of $X$ is $n-4$.

In this paper, we consider the finite morphisms from abelian surfaces to nonsingular complete toric $4$-folds satisfying the following condition.

\begin{Def}\label{defoftndabel}
{\rm
Let $X$ be a $4$-dimensional complete nonsingular toric variety and $A$ an abelian surface. A finite morphism $\varphi:A\rightarrow X$ is called a {\em totally nondegenerate finite morphism} if $D_i\cap \varphi(A)$ is non-empty on $\varphi(A)$ for any $T_N$-invariant prime divisor $D_i$ $(1\leq i\leq n)$. If $\varphi$ is an embedding, we call $\varphi$ a {\em totally nondegenerate embedding}.
}

\end{Def}

\begin{Not}

{\rm
We use the following notation throughout this paper.
\begin{enumerate}
\item Let $\varphi:A\rightarrow X$ be a totally nondegenerate finite morphism. For a $T_N$-invariant prime divisor $D_i$ on $X$, we put $C_i:=\varphi^{\ast}D_i$. $C_i$ is a divisor on $A$.
\item The types of nonsingular toric Fano $4$-folds are in the sense of Batyrev \cite{batyrev4} and Sato \cite{sato1}. We use characters $\mathcal{B},\mathcal{C},\mathcal{D}\ldots$ instead of $B,C,D\ldots$. See Table $1$ in Section \ref{tablefano}.
\end{enumerate}
}
\end{Not}

For a basis $\{x_1,x_2,x_3,x_4\}$ of $N$, by computing the divisors of the rational functions ${\bf e}(x^{\ast}_1),$ ${\bf e}(x^{\ast}_2),$ ${\bf e}(x^{\ast}_3),$ ${\bf e}(x^{\ast}_4)\in {\bf C}(X)$, where $\{x^{\ast}_1,x^{\ast}_2,x^{\ast}_3,x^{\ast}_4\}\subset M$ is the dual basis of $\{x_1,x_2,x_3,x_4\}$, we get four linear relations
$${\rm div}\left({\bf e}(x^{\ast}_1)\right)=0,\ {\rm div}\left({\bf e}(x^{\ast}_2)\right)=0,\ {\rm div}\left({\bf e}(x^{\ast}_3)\right)=0\mbox{ and }{\rm div}\left({\bf e}(x^{\ast}_4)\right)=0$$
among $T_N$-invariant prime divisors in ${\rm Pic}(X)$. We often use this argument in the following sections.

\section{Criteria for non-existence}\label{somecri}

\hspace{5mm} In this section, we present criteria for the non-existence of totally nondegenerate finite morphism from abelian surface $A$ to a projective nonsingular toric 4-fold. We reduce Kajiwara's method in \cite{kajiwara2} and \cite{kajiwara1} to a more convenient form. For fundamental properties of primitive collections and primitive relations, see Batyrev \cite{batyrev3}, \cite{batyrev4} and Sato \cite{sato1}.

\begin{Lem}\label{effonabel}
Let $A$ be an abelian surface and $D$ an effective divisor on $A$. Then we have $D^2\geq 0$.

\end{Lem}

\proof
We may assume that $D$ is an irreducible curve on $A$. For some point $x\in A$, we have $D(D+x)\geq 0$, where $D+x$ is the translation of $D$ by $x$. Since $D$ and $D+x$ are algebraically equivalent, we have $D^2=D(D+x)\geq 0$.\qed

\begin{Lem}\label{conntocomp}
Let $X$ be a complete nonsingular toric $4$-fold and $\varphi:A\rightarrow X$ a totally nondegenerate finite morphism. If $(\varphi^{\ast}D_i)(\varphi^{\ast}D_j)=0$ and $(\varphi^{\ast}D_j)(\varphi^{\ast}D_k)=0$ where $1\leq i,j,k\leq n$, then we have $(\varphi^{\ast}D_i)(\varphi^{\ast}D_k)=0$.

\end{Lem}

\proof
Suppose that $(\varphi^{\ast}D_i)(\varphi^{\ast}D_k)>0$. Since $\left( \varphi^{\ast}D_i+\varphi^{\ast}D_k\right)^{2}=(\varphi^{\ast}D_i)^{2}+(\varphi^{\ast}D_k)^{2}+2(\varphi^{\ast}D_i)(\varphi^{\ast}D_k)\geq 2(\varphi^{\ast}D_i)(\varphi^{\ast}D_k)>0$ by Lemma \ref{effonabel}, $\varphi^{\ast}D_i+\varphi^{\ast}D_k$ is an ample divisor on $A$. On the other hand, $\left( \varphi^{\ast}D_i+\varphi^{\ast}D_k\right)(\varphi^{\ast}D_j)=0$, by assumption. This contradicts the fact that $\varphi^{\ast}D_i+\varphi^{\ast}D_k$ is ample. Therefore, $(\varphi^{\ast}D_i)(\varphi^{\ast}D_k)=0$.\hfill q.e.d.

\bigskip

For a totally nondegenerate finite morphism $\varphi:A\rightarrow X$, we define a graph $\Gamma_{\varphi}$ as follows: The vertex set of $\Gamma_{\varphi}$ is $\{ 1,\ldots,n\}$, and $\{ i,j\}$ is an edge of $\Gamma_{\varphi}$ if $i\neq j$ and $(\varphi^{\ast}D_i)(\varphi^{\ast}D_j)=C_iC_j=0$.

\begin{Rem}\label{conntocomprem}
{\rm
Lemma \ref{conntocomp} implies that every connected component of $\Gamma_{\varphi}$ is complete, i.e., any pair of distinct vertices is connected by an edge. Especially, if $\Gamma_{\varphi}$ is connected, then $\Gamma_{\varphi}$ is a complete graph.
}

\end{Rem}

\begin{Lem}\label{projtononcomp}
Let $X$ be a {\em projective} nonsingular toric $4$-fold. If there exists a totally nondegenerate finite morphism $\varphi:A\rightarrow X$, then $(\varphi^{\ast}D_i)(\varphi^{\ast}D_j)>0$ for some $1\leq i<j\leq n$.

\end{Lem}

\proof
Since $X$ is projective, there exists an ample effective divisor $\sum_{k=1}^{n}a_kD_k$ on $X$. Since $\sum_{k=1}^{n}a_k\varphi^{\ast}D_k$ is also ample on $A$, we have $\left( \sum_{k=1}^{n}a_k\varphi^{\ast}D_k\right)^{2}>0$. Therefore, we have $(\varphi^{\ast}D_i)(\varphi^{\ast}D_j)>0$ for some $1\leq i\leq j\leq n$. If $i=j$, then there exists $1\leq l\leq n$ such that $l\neq i$ and $(\varphi^{\ast}D_i)(\varphi^{\ast}D_l)>0$ by Lemma \ref{conntocomp}.\hfill q.e.d.

\bigskip

\begin{Rem}\label{keypoint}
{\rm
Lemma \ref{projtononcomp} implies that $\Gamma_{\varphi}$ is not complete. Especially, by Remark \ref{conntocomprem}, $\Gamma_{\varphi}$ is not connected.
}

\end{Rem}

\begin{Rem}\label{keypoint2}
{\rm
For $n>5$, the assertion in Remark \ref{keypoint} is also true if we replace the vertex set of $\Gamma_{\varphi}$ by $S\subset\{1,\ldots,n\}$ such that $\left\{ D_i\right\}_{i\in S} \subset {\rm Pic}(X)$ generates ${\rm Pic}(X)$.
}

\end{Rem}

By using this incompleteness of $\Gamma_{\varphi}$, we can show the non-existence of totally nondegenerate finite morphisms for some projective nonsingular toric 4-folds. For example, the following holds.

\begin{Ex}\label{g1noexist}
{\rm
Let $X$ be the nonsingular projective toric $4$-fold corresponding to the fan $\Sigma$ whose primitive relations are
$$x_1+x_7=0,\ x_2+x_3+x_4=ax_1,\ x_4+x_5+x_6=(a+1)x_1,$$
$$x_5+x_6+x_7=x_2+x_3\ \mbox{and}\ x_1+x_2+x_3=x_5+x_6,$$
where $\G(\Sigma)=\{x_1,\ldots,x_7\}$ and $a$ is a positive integer. $D_1D_7=0$ on $X$ because $\{x_1,x_7\}$ is a primitive collection of $\Sigma$, and by a basis $\{x_1,x_2,x_4,x_5\}$ of $N$, we have
$$(1)\ D_1+aD_3+(a+1)D_6-D_7=0,\ (2)\ D_2-D_3=0,$$
$$(3)\ -D_3+D_4-D_6=0\ \mbox{and}\ (4)\ D_5-D_6=0$$
in ${\rm Pic}(X)$, respectively. Suppose that there exists a totally nondegenerate finite morphism $\varphi:A\rightarrow X$. By intersecting $D_1$ with both sides of $(1)$ and restricting the result to $A$, we have $C_1^{2}+aC_1C_3+(a+1)C_1C_6-C_1C_7=0$. So $C_1C_3=C_1C_6=0$, and $C_1C_2=C_1C_4=C_1C_5=0$ by $(2)$, $(3)$ and $(4)$. Hence $\Gamma_{\varphi}$ is connected, a contradiction to Remark \ref{keypoint}. Therefore, $X$ admits no totally nondegenerate finite morphism.
}

\end{Ex}

\begin{Rem}
{\rm
In Example \ref{g1noexist}, if $a=1$, then $X$ is the nonsingular toric Fano $4$-fold of type $\mathcal{G}_1$. So there exists no totally nondegenerate finite morphism to the nonsingular toric Fano $4$-fold of type $\mathcal{G}_1$.
}

\end{Rem}

\begin{Prop}\label{ainfdivcont}
If a nonsingular projective toric $4$-fold $X$ has an equivariant projective birational divisorial contraction to a possibly singular point, then there exist no totally nondegenerate finite morphism to $X$.

\end{Prop}

\proof
Suppose that there exists a totally nondegenerate finite morphsim $\varphi:A\rightarrow X$. By the Mori theory for projective toric varieties (see Reid \cite{reid1}), we may assume, without loss of generality, that we have a primitive relation $x_1+x_2+x_3+x_4=ax_5$, where $a$ is a positive integer. Obviously $D_5D_i=0$ for $6\leq i\leq n$. By a basis $\{x_1,x_2,x_3,x_5\}$ of $N$, we have $aD_4+D_5+b_6D_6+\cdots+b_nD_n=0$ in ${\rm Pic}(X)$. So $aC_4C_5+C_5^{2}=0$ on $A$. Therefore $C_4C_5=0$, and similarly $C_iC_5=0$ for $1\leq i\leq 4$. This means that $\Gamma_{\varphi}$ is connected, a contradiction to Remark \ref{keypoint}.\hfill q.e.d.

\begin{Rem}
{\rm
By Proposition \ref{ainfdivcont}, there exist no totally nondegenerate finite morphism to the nonsingular toric Fano $4$-folds of types $\mathcal{B}_1$, $\mathcal{B}_2$, $\mathcal{B}_3$ and $\mathcal{E}_3$.
}

\end{Rem}

We now consider the case where $X$ is decomposed into the product of ${\bf P}^{1}$ and a projective nonsingular toric $3$-fold. Let $X={\bf P}^{1}\times X'$ be a projective nonsingular toric $4$-fold, where $X'$ is a projective nonsingular toric $3$-fold. Suppose that $x_1$ and $x_2$ correspond to the class of fibers of the first projection $X\rightarrow {\bf P}^{1}$, where $\G(\Sigma)=\{ x_1,x_2,\ldots,x_n\}$. Then the following holds.

\begin{Prop}\label{subconn}
Suppose that there exists a totally nondegenerate finite morphism $\varphi:A\rightarrow X$. We define a subgraph $\Gamma'$ of $\Gamma_{\varphi}$ as follows$:$ The vertex set of $\Gamma'$ is $\{ 3,\ldots,n\}$, and $\{i,j\}$ $(3\leq i<j\leq n)$ is an edge of $\Gamma'$ if $\{i,j\}$ is an edge of $\Gamma_{\varphi}$. Then $\Gamma'$ is not complete.

\end{Prop}

\proof
Since each fiber of the second projection $p_2:X\rightarrow X'$ is ${\bf P}^1$, each fiber of $p_2$ is not contained in the abelian surface $A$. So $p_2\circ\varphi$ is a finite morphism. Therefore, for an ample divisor $E=\sum_{i=3}^{n}a_iE_i$ on the projective variety $X'$, where $E_3,\ldots,E_n$ are the toric divisors corresponding to $x_3,\ldots,x_n$, respectively, $(p_2\circ\varphi)^{\ast}(E)=\sum_{i=3}^{n}a_i(p_2\circ\varphi)^{\ast}D_i$ is also an ample divisor on $A$. So there exist $3\leq i<j\leq n$ such that $((p_2\circ\varphi)^{\ast}D_i)((p_2\circ\varphi)^{\ast}D_j)\neq 0$. Since $\{i,j\}$ is not a edge of $\Gamma'$, the graph $\Gamma'$ is not complete.\hfill q.e.d.

\bigskip

\begin{Ex}\label{l5noexist}
{\rm
Let $X$ be the nonsingular projective toric $4$-fold corresponding to the fan $\Sigma$ whose primitive relations are
$$x_1+x_8=0,\ x_2+x_3=0,\ x_4+x_5=ax_3\ \mbox{and}\ x_6+x_7=ax_3,$$
where $\G(\Sigma)=\{x_1,\ldots,x_8\}$ and $a$ is a positive integer. $D_1D_8=D_2D_3=D_4D_5=D_6D_7=0$ on $X$, and by a basis $\{x_1,x_2,x_4,x_6\}$ of $N$, we have
$$(1)\ D_1-D_8=0,\ (2)\ D_2-D_3-aD_5-aD_7=0,\ (3)\ D_4-D_5=0\ \mbox{and}\ (4)\ D_6-D_7=0$$
in ${\rm Pic}(X)$, respectively. $X$ is isomorphic to ${\bf P}^{1}\times X'$, where $X'$ is a toric $3$-fold, and $D_1$ and $D_8$ are fibers of the first projection $X\rightarrow {\bf P}^{1}$. Suppose that there exists a totally nondegenerate finite morphism $\varphi:A\rightarrow X$. By $(2)$, we have $C_3^{2}+aC_3C_5+aC_3C_7=C_2C_3=0$, hence $C_3C_5=C_3C_7=0$. Consequently, $C_3C_4=C_3C_6=0$ by $(3)$ and $(4)$. Hence the graph $\Gamma'$ as in Proposition \ref{subconn} is connected, a contradiction to  Proposition \ref{subconn}. Therefore, $X$ admits no totally nondegenerate finite morphism.
}

\end{Ex}

\begin{Rem}
{\rm
In Example \ref{l5noexist}, if $a=1$, then $X$ is the nonsingular toric Fano $4$-fold of type $\mathcal{L}_5$. So there exists no totally nondegenerate finite morphism to the nonsingular toric Fano $4$-fold of type $\mathcal{L}_5$.
}

\end{Rem}

For the main theorem of this paper, we show some results for the non-existence of totally nondegenerate finite morphisms, using Remark \ref{keypoint} and Proposition \ref{subconn}.

\begin{Prop}\label{faoverp2}
Let $X$ be an $F_a$-bundle over ${\bf P}^{2}$, where ${\rm F}_a$ is the Hirzebruch surface of degree $a$ $(a\geq 0)$, and $\G(\Sigma)=\{ x_1,\ldots x_7\}$. We introduce a coordinate in $N$ so that the coordinates of $x_1,x_2,x_3,x_4,x_5,x_6$ and $x_7$ are
\[
\pmatrix{
1 \cr
0 \cr
0 \cr
0 \cr},\ 
\pmatrix{
0 \cr
1 \cr
0 \cr
0 \cr},\ 
\pmatrix{
-1 \cr
-1 \cr
s \cr
t \cr},\ 
\pmatrix{
0 \cr
0 \cr
1 \cr
0 \cr},\ 
\pmatrix{
0 \cr
0 \cr
-1 \cr
a \cr},\ 
\pmatrix{
0 \cr
0 \cr
0 \cr
1 \cr}\mbox{ and }
\pmatrix{
0 \cr
0 \cr
0 \cr
-1 \cr},
\]
respectively, where $s$ and $t$ are integers. In this situation, the following hold:
\begin{enumerate}
\item If $s=t=0$, then $X$ is isomorphic to ${\bf P}^{2}\times{\rm F}_a$.
\item In the case $s\neq 0$ or $t\neq 0$, if one of the following conditions is satisfied, then $X$ admits no totally nondegenerate finite morphism.
\begin{enumerate}
\item $a=0$.
\item $a>0$ and $t\geq 0$.
\item $a>0$, $s>0$, $t<0$ and $as+t\geq 0$.
\end{enumerate}
\end{enumerate}

\end{Prop}

\proof
Suppose that there exists a totally nondegenerate finite morphism $\varphi:A\rightarrow X$. $($i$)$ is trivial, so let $s\neq 0$ or $t\neq 0$. By a basis $\{x_1,x_2,x_4,x_6\}$ of $N$, we have
$$(1)\ D_1-D_3=0,\ (2)\ D_2-D_3=0,\ (3)\ sD_3+D_4-D_5=0\mbox{ and }(4)\ tD_3+aD_5+D_6-D_7=0$$
in ${\rm Pic}(X)$, respectively. Moreover, we have $D_4D_5=D_6D_7=0$ on $X$, and $C_4^{2}=C_5^{2}=C_6^{2}=C_7^{2}=0$ on $A$.

(a) Let $a=0$. In this case, $X$ is isomorphic to ${\bf P}^{1}\times X'$, where $X'$ is a toric $3$-fold.

If $s=0$ and $t\neq 0$, then $D_4$ and $D_5$ are in the class of fibers of the first projection $X\rightarrow{\bf P}^{1}$. Since $tD_3+D_6-D_7=0$ by $(4)$, we have $tC_3C_6=-C_6^{2}+C_6C_7=0$. Therefore $C_3C_6=0$, and $C_1C_6=C_2C_6=0$ by $(1)$ and $(2)$. This contradicts Proposition \ref{subconn}.

If $s\neq 0$ and $t=0$, then $D_6$ and $D_7$ are in the class of fibers of the first projection $X\rightarrow{\bf P}^{1}$. Since $sC_3C_4=-C_4^{2}+C_4C_5=0$ by $(3)$, we have $C_3C_4=0$. On the other hand, $C_1C_4=C_2C_4=0$ by $(1)$ and $(2)$. This contradicts Proposition \ref{subconn}.

(b) Let $a>0$ and $t\geq 0$. Since
$$(5)\ tC_3C_6+aC_5C_6=-C_6^{2}+C_6C_7=0$$
by $(4)$, we have $C_5C_6=0$.

If $t>0$, then $C_3C_6=0$ by $(5)$. Moreover $C_1C_6=C_2C_6=0$ by $(1)$ and $(2)$. So $\Gamma_{\varphi}$ is connected, a contradiction to Remark \ref{keypoint}.

Let $t=0$. Then $s\neq 0$ by assumption. So we have $C_3C_4=0$ as above, and $C_1C_4=C_2C_4=0$ by $(1)$ and $(2)$. So $\Gamma_{\varphi}$ is connected, a contradiction to Remark \ref{keypoint}.

(c) Let $a>0$, $s>0$, $t<0$ and $as+t\geq 0$. Then we have $C_3C_4=C_1C_4=C_2C_4=0$ as above. On the other hand, by $(3)$ and $(4)$, we have
$$(6)\ -tD_4+(as+t)D_5+sD_6-sD_7=0$$
in ${\rm Pic}(X)$. So $-tC_4C_6+(as+t)C_5C_6=-sC_6^{2}+sC_6C_7=0$, and we have $C_4C_6=0$ by the assumptions $t<0$ and $as+t\geq 0$. Therefore $\Gamma_{\varphi}$ is connected, a contradiction to Remark \ref{keypoint}.\hfill q.e.d.

\begin{Rem}
{\rm
By Proposition \ref{faoverp2}, there exists no totally nondegenerate finite morphism to the nonsingular toric Fano $4$-folds of types $\mathcal{D}_1$, $\mathcal{D}_2$, $\mathcal{D}_3$, $\mathcal{D}_5$, $\mathcal{D}_6$, $\mathcal{D}_8$, $\mathcal{D}_9$, $\mathcal{D}_{12}$ and $\mathcal{D}_{16}$. The corresponding $a$, $s$ and $t$ are as follows:

\begin{center}
\begin{tabular}{|c||c|c|c|c|c|c|c|c|c|}\hline
 & $\mathcal{D}_1$ & $\mathcal{D}_2$ & $\mathcal{D}_3$ & $\mathcal{D}_5$ & $\mathcal{D}_6$ & $\mathcal{D}_8$ & $\mathcal{D}_9$ & $\mathcal{D}_{12}$ & $\mathcal{D}_{16}$ \\ \hline
$a$ & 1 & 1 & 1 & 0 & 1 & 1 & 0 & 0 & 1 \\ \hline
$s$ & 0 & 2 & 1 & 2 & 0 & 1 & 1 & 1 & 1 \\ \hline
$t$ & 2 & 0 & 1 & 0 & 1 & 0 & 1 & 0 & -1 \\ \hline
\end{tabular}
\end{center}

}

\end{Rem}

\section{2-blow-up}\label{main2blowup}

\hspace{.5cm} The following is useful for deriving the main result in Section $6$.

\begin{Prop}\label{2blowupprop}
Let $X$ and $\widetilde{X}$ be nonsingular projective toric $4$-folds and $\psi:\widetilde{X}\rightarrow X$ a $2$-blow-up, where a ``$2$-blow-up'' means an equivariant blow-up along a $T_N$-invariant subvariety of codimension $2$. If $X$ admits no totally nondegenerate finite morphism, then $\widetilde{X}$ admits no totally nondegenerate finite morphism either.

\end{Prop}

\proof
If there exists a totally nondegenerate finite morphism $\varphi:A\rightarrow \widetilde{X}$, then $\psi\circ\varphi:A\rightarrow X$ is also a totally nondegenerate finite morphism (see Mumford \cite{mumford1}, p. 88). This is a contradiction.\qed

\bigskip

Especially, we have the following.

\begin{Cor}\label{dia2}
Let $X_1\leftarrow X_2\leftarrow\cdots\leftarrow X_{n-1}\leftarrow X_n$ be a sequence of $2$-blow-ups among nonsingular projective toric $4$-folds. If $X_1$ admits no totally nondegenerate finite morphism then $X_n$ admits no totally nondegenerate finite morphism either.

\end{Cor}

We end this section by proposing the following conjecture.

\begin{Conj}\label{2blowupconj}
Let $X$ and $\widetilde{X}$ be nonsingular projective toric $4$-folds and $\psi:\widetilde{X}\rightarrow X$ a $2$-blow-up. If $X$ admits no totally nondegenerate embedding, then $\widetilde{X}$ admits no totally nondegenerate embedding either.

\end{Conj}

\section{Some examples}\label{secex}

\hspace{.5cm} In this section, to describe the main result in Section \ref{mainresult}, we show the non-existence of totally nondegenerate embeddings for some other nonsingular toric Fano $4$-folds.

\bigskip

(a) ``type $\mathcal{I}$'s''\quad Let $X$ be the nonsingular projective toric $4$-fold corresponding to the fan $\Sigma$ defined as follows: Let $\G(\Sigma)=\{x_1,\ldots,x_8\}\subset N$ such that the coordinates of $x_1,\ldots,x_8$ are
\[
\pmatrix{
1 \cr
0 \cr
0 \cr
0 \cr},\ 
\pmatrix{
-1 \cr
b \cr
0 \cr
c \cr},\ 
\pmatrix{
0 \cr
1 \cr
0 \cr
0 \cr},\ 
\pmatrix{
0 \cr
0 \cr
1 \cr
0 \cr},\ 
\pmatrix{
0 \cr
-1 \cr
-1 \cr
a+1 \cr},\ 
\pmatrix{
0 \cr
0 \cr
0 \cr
-1 \cr},\ 
\pmatrix{
0 \cr
1 \cr
0 \cr
-1 \cr},\ 
\pmatrix{
0 \cr
0 \cr
0 \cr
1 \cr},
\]
respectively, and that the primitive collections of $\Sigma$ are $\{ x_3,x_4,x_5\}$, $\{ x_4,x_5,x_7\}$, $\{ x_7,x_8\}$, $\{ x_3,x_6\}$, $\{ x_6,x_8\}$ and $\{ x_1,x_2\}$. For some values of $a$, $b$ and $c$, $X$ becomes the nonsingular toric Fano $4$-fold of type $\mathcal{I}$. The corresponding $a$, $b$ and $c$ are as follows:
\[
\begin{tabular}{|c||c|c|c|c|}\hline
 & $\mathcal{I}_4$ & $\mathcal{I}_6$ & $\mathcal{I}_{12}$ & $\mathcal{I}_{15}$ \\ \hline
$a$ & 1 & 0 & 0 & 1 \\ \hline
$b$ & 1 & 1 & 0 & 0 \\ \hline
$c$ & -1 & 0 & -1 & -1 \\ \hline
\end{tabular}
\]
$D_7D_8=D_3D_6=D_6D_8=D_1D_2=0$ on $X$, and by a basis $\{x_1,x_3,x_4,x_8\}$ of $N$, we have
$$(1)\ D_1-D_2=0,\ (2)\ bD_2+D_3-D_5+D_7=0,$$
$$(3)\ D_4-D_5=0\ \mbox{and}\ (4)\ cD_2+(a+1)D_5-D_6-D_7+D_8=0$$
in ${\rm Pic}(X)$, respectively. Moreover, we have
$$(5)\ (b+c)D_2+D_3+aD_5-D_6+D_8=0$$
by $(2)$ and $(4)$. Suppose that there exists a totally nondegenerate finite morphism $\varphi:A\rightarrow X$.

If $b=0$ and $c=-1$, then $C_5C_6=C_3C_6+C_6C_7=0$ by $(2)$. On the other hand, by $(4)$, we have $C_2C_5=C_5C_8+(a+1)C_5^{2}-C_5C_6-C_5C_7=0$. Hence $\Gamma_{\varphi}$ is connected, a contradiction to Remark \ref{keypoint}. Therefore, $X$ admits no totally nondegenerate finite morphism. Especially, the nonsingular toric Fano $4$-folds of types $\mathcal{I}_{12}$ and $\mathcal{I}_{15}$ admit no totally nondegenerate finite morphism.

In the case $b=1$, if $a=1$ and $b+c=0$, then $C_3C_5=-C_3^{2}+C_3C_6-C_3C_8=0$ by $(5)$, and $C_2C_3=-C_3^{2}+C_3C_5$ $-C_3C_7=0$ by $(2)$. On the other hand, if $a=0$ and $b+c=1$, then $C_2C_3=0$ by $(5)$, and $C_3C_5=0$ by $(2)$. In any case, $\Gamma_{\varphi}$ is connected, a contradiction to Remark \ref{keypoint}. Therefore, $X$ admits no totally nondegenerate finite morphism. Especially, the nonsingular toric Fano $4$-folds of types $\mathcal{I}_{4}$ and $\mathcal{I}_{6}$ admit no totally nondegenerate finite morphism.

\bigskip

(b) ``type $\mathcal{J}_2$''\quad Let $X$ be the nonsingular projective toric $4$-fold corresponding to the fan $\Sigma$ whose primitive relations are
$$x_3+x_6=x_7,\ x_1+x_2+x_8=x_4+x_5,\ x_4+x_5+x_6=x_1+x_2,\ x_7+x_8=x_3,\ x_6+x_8=0,$$
$$x_3+x_4+x_5=x_8,\ x_4+x_5+x_7=0,\ x_1+x_2+x_3=0\ \mbox{and}\ x_1+x_2+x_7=x_6,$$
where $\G(\Sigma)=\{x_1,\ldots,x_8\}$. $X$ is the nonsingular toric Fano $4$-fold of type $\mathcal{J}_2$. $D_3D_6=D_7D_8=D_6D_8=0$ on $X$, and by a basis $\{x_1,x_2,x_4,x_5\}$ of $N$, we have
$$(1)\ D_1-D_3+D_6-D_8=0,\ (2)\ D_2-D_3+D_6-D_8=0,$$
$$(3)\ D_4-D_6-D_7+D_8=0\ \mbox{and}\ (4)\ D_5-D_6-D_7+D_8=0$$
in ${\rm Pic}(X)$, respectively. So we have $D_1=D_2$ and $D_4=D_5$. Suppose that there exists a totally nondegenerate finite morphism $\varphi:A\rightarrow X$. By $(1)$, we have $C_1C_3=C_3^{2}-C_3C_6+C_3C_8=0$. On the other hand, by $(3)$, we have $C_3C_4=C_3C_6+C_3C_7-C_3C_8=0$. Hence $\Gamma_{\varphi}$ is connected, a contradiction to Remark \ref{keypoint}. Therefore, $X$ admits no totally nondegenerate finite morphism.

\bigskip

(c) ``type $\mathcal{L}$'s''\quad Let $X$ be the nonsingular projective toric $4$-fold corresponding to the fan $\Sigma$ defined as follows: Let $\G(\Sigma)=\{x_1,\ldots,x_8\}\subset N$ such that the coordinates of $x_1,\ldots,x_8$ are
\[
\pmatrix{
0 \cr
1 \cr
0 \cr
0 \cr},\ 
\pmatrix{
-1 \cr
1 \cr
0 \cr
0 \cr},\ 
\pmatrix{
1 \cr
0 \cr
0 \cr
0 \cr},\ 
\pmatrix{
0 \cr
0 \cr
1 \cr
0 \cr},\ 
\pmatrix{
a \cr
b \cr
-1 \cr
0 \cr},\ 
\pmatrix{
0 \cr
0 \cr
0 \cr
1 \cr},\ 
\pmatrix{
c \cr
d \cr
0 \cr
-1 \cr},\ 
\pmatrix{
0 \cr
-1 \cr
0 \cr
0 \cr},
\]
respectively, and that the primitive collections of $\Sigma$ are $\{ x_1,x_8\}$, $\{ x_2,x_3\}$, $\{ x_4,x_5\}$ and $\{ x_6,x_7\}$. For some values of $a$, $b$, $c$ and $d$, $X$ becomes the nonsingular toric Fano $4$-fold of type $\mathcal{L}$. The corresponding $a$, $b$, $c$ and $d$ are as follows:
\[
\begin{tabular}{|c||c|c|c|}\hline
 & $\mathcal{L}_1$ & $\mathcal{L}_{2}$ & $\mathcal{L}_{10}$ \\ \hline
$a$ & 0 & 1 & 1 \\ \hline
$b$ & 1 & 0 & 0 \\ \hline
$c$ & 0 & 1 & -1 \\ \hline
$d$ & 1 & 0 & 1 \\ \hline
\end{tabular}
\]
$D_1D_8=D_2D_3=D_4D_5=D_6D_7=0$ on $X$, and by a basis $\{x_1,x_3,x_4,x_6\}$ of $N$, we have
$$(1)\ D_1+D_2+bD_5+dD_7-D_8=0,\ (2)\ -D_2+D_3+aD_5+cD_7=0,$$
$$(3)\ D_4-D_5=0\ \mbox{and}\ (4)\ D_6-D_7=0$$
in ${\rm Pic}(X)$, respectively. Suppose that there exists a totally nondegenerate finite morphism $\varphi:A\rightarrow X$. By $(1)$ and $(2)$, we have $C_1C_2+bC_1C_5+dC_1C_7=-C_1^{2}+C_1C_8=0$ and $aC_3C_5+cC_3C_7=-C_3^{2}+C_2C_3=0$, respectively. Especially, $C_1C_2=0$.

If either $a>0$ and $c>0$ or $b>0$ and $d>0$, then $C_3C_5=C_3C_5=0$ or $C_1C_5=C_1C_7=0$, respectively. In any case, $\Gamma_{\varphi}$ is connected, a contradiction to Remark \ref{keypoint}. Therefore, $X$ admits no totally nondegenerate finite morphism. Especially, the nonsingular toric Fano $4$-folds of types $\mathcal{L}_{1}$ and $\mathcal{L}_{2}$ admit no totally nondegenerate finite morphism.

If $X$ is of type $\mathcal{L}_{10}$, then $C_1C_7=0$ because $b\geq 0$ and $d>1$. Moreover, by $(2)$, we have $C_1C_5=-C_1C_3+C_2C_3+C_1C_7=0$. Hence $\Gamma_{\varphi}$ is connected, a contradiction to Remark \ref{keypoint}. Therefore, the nonsingular toric Fano $4$-fold of type $\mathcal{L}_{10}$ admits no totally nondegenerate finite morphism.

\bigskip

(d) ``type $\mathcal{L}_{12}$''\quad This case is special. Let $X$ be the nonsingular projective toric $4$-fold corresponding to the fan $\Sigma$ whose primitive relations are
$$x_1+x_8=0,\ x_2+x_3=x_1,\ x_4+x_5=x_8\ \mbox{and}\ x_6+x_7=x_4,$$
where $\G(\Sigma)=\{x_1,\ldots,x_8\}$. $X$ is the nonsingular toric Fano $4$-fold of type $\mathcal{L}_{12}$. $D_1D_8=D_2D_3=D_4D_5=D_6D_7=0$ on $X$, and by a basis $\{x_1,x_2,x_4,x_6\}$ of $N$, we have
$$(1)\ D_1+D_3-D_5-D_8=0,\ (2)\ D_2-D_3=0,$$
$$(3)\ D_4-D_5+D_7=0\ \mbox{and}\ (4)\ D_6-D_7=0$$
in ${\rm Pic}(X)$, respectively. So we have $D^{2}_2=D^{2}_3=D^{2}_6=D^{2}_7=0$,
$$(5)\ D^{2}_{8}=D_3D_8-D_5D_8\ \mbox{and}\ (6)\ D^{2}_5=D_5D_7$$
on $X$. Suppose that there exists a totally nondegenerate embedding $\varphi:A\hookrightarrow X$. Then by $(3)$, we have $D_4D_7A=-(D_4A)^{2}+D_4D_5A=0$. Moreover, since $D_1A$ is an effective divisor on $A$ and $D_3D_8A-D_5D_8A=0$ by $(5)$, we have $(D_1A)^{2}=(-D_3A+D_5A+D_8A)^{2}=-2D_3D_5A-2D_3D_8A+2D_5D_8A=-2D_3D_5A\geq 0$. Therefore, we have $D_3D_5A=0$. On the other hand, $\{ D_3D_5,D_3D_7,D_3D_8,D_5D_7,D_5D_8,D_7D_8\}$ generates ${\rm A}^{2}(X)$ by the equalities $(5)$ and $(6)$. So we can express the class of $A$ in ${\rm A}^{2}(X)$ as
$$A=a_1D_3D_5+a_2D_3D_7+a_3D_3D_8+a_4D_5D_7+a_5D_5D_8+a_6D_7D_8\in {\rm A}^{2}(X).$$
Since $D_3D^{2}_5D_7=D_3(D_4+D_7)D_5D_7=0$, $D_3D^{2}_5D_8=D_3(D_4+D_7)D_5D_8=1$ and $D^{2}_5D_7D_8=(D_4+D_7)D_5D_7D_8=0$, we have the following:
\[
(7)\ \left\{
\begin{array}{c}
D_3D_5A=a_4D_3D^{2}_5D_7+a_5D_3D^{2}_5D_8+a_6D_3D_5D_7D_8=a_5+a_6=0, \\
D_3D_7A=a_5D_3D_5D_7D_8=a_5=0\mbox{ and} \\
D_5D_7A=a_1D_3D^{2}_5D_7+a_3D_3D_5D_7D_8+a_5D^{2}_5D_7D_8=a_3=0. \\
\end{array}
\right.
\]
By these equalities, $a_3=a_5=a_6=0$. So $A^{2}=a_1D_3D_5A+a_2D_3D_7A+a_4D_5D_7A=0$. Therefore, by the following, we have $c_2(X)A=0$.

\begin{Lem}[Van de Ven \cite{vandeven1}, Proposition 3]\label{selfintersection}
Let $A\hookrightarrow X$ be an embedding from an abelian surface to a $4$-dimensional nonsingular projective toric variety. Then we have
$c_2(X)A=A^{2}$ in the group ${\rm A}^{4}(X)$ of codimension four cycles on $X$ modulo rational equivalence.

\end{Lem}

Since
$$c_2(X)=\sum_{1\leq i<j\leq n}D_iD_j,$$
$\Gamma_{\varphi}$ is connected. This contradicts Remark \ref{keypoint}. So $X$ admits no totally nondegenerate embedding.

\bigskip

(e) ``type $\mathcal{M}$'s''\quad Let $X$ be the nonsingular projective toric $4$-fold corresponding to the fan $\Sigma$ defined as follows: Let $\G(\Sigma)=\{x_1,\ldots,x_7\}\subset N$ such that the coordinates of $x_1,\ldots,x_7$ are
\[
\pmatrix{
1 \cr
0 \cr
0 \cr
0 \cr},\ 
\pmatrix{
0 \cr
1 \cr
0 \cr
0 \cr},\ 
\pmatrix{
-1 \cr
-1 \cr
1 \cr
1 \cr},\ 
\pmatrix{
0 \cr
0 \cr
1 \cr
0 \cr},\ 
\pmatrix{
a \cr
0 \cr
-1 \cr
0 \cr},\ 
\pmatrix{
0 \cr
0 \cr
0 \cr
1 \cr},\ 
\pmatrix{
ac+b \cr
0 \cr
-c \cr
-1 \cr},\ 
\pmatrix{
-1 \cr
0 \cr
0 \cr
0 \cr},
\]
respectively, and that the primitive collections of $\Sigma$ are $\{ x_1,x_8\}$, $\{ x_1,x_2,x_3\}$, $\{ x_4,x_6,x_8\}$, $\{ x_4,x_5\}$, $\{ x_6,x_7\}$, $\{ x_2,x_3,x_5\}$ and $\{ x_2,x_3,x_7\}$. For some values of $a$, $b$ and $c$, $X$ becomes the nonsingular toric Fano $4$-fold of type $\mathcal{M}$. The corresponding $a$, $b$ and $c$ are as follows:
\[
\begin{tabular}{|c||c|c|c|c|}\hline
 & $\mathcal{M}_1$ & $\mathcal{M}_2$ & $\mathcal{M}_3$ & $\mathcal{M}_4$ \\ \hline
$a$ & 0 & 1 & 1 & 1 \\ \hline
$b$ & 0 & 1 & 0 & 0 \\ \hline
$c$ & 0 & 0 & 1 & 0 \\ \hline
\end{tabular}
\]
$D_1D_8=D_4D_5=D_6D_7=0$ on $X$, and by a basis $\{x_1,x_2,x_4,x_6\}$ of $N$, we have
$$(1)\ D_1-D_3+aD_5+(ac+b)D_7-D_8=0,\ (2)\ D_2-D_3=0,$$
$$(3)\ D_3+D_4-D_5-cD_7=0\ \mbox{and}\ (4)\ D_3+D_6-D_7=0$$
in ${\rm Pic}(X)$, respectively. Suppose that there exists a totally nondegenerate finite morphism $\varphi:A\rightarrow X$. By $(4)$, we have $C_3C_6=-C_6^{2}+C_6C_7=0$. By $(1)$, $(3)$ and $(4)$, we have
$$(5)\ D_1+(a-1)D_3+aD_4+bD_7-D_8=0\ \mbox{and}\ (6)\ D_1+aD_5+D_6+(ac+b-1)D_7-D_8=0.$$
Let $a\geq 0$, $b\geq 0$ and $c\geq 0$.

If $a=1$ and $c=0$, then by $(3)$, we have $C_3C_4=-C_4^{2}+C_4C_5=0$. On the other hand, by $(5)$, we have $C_1C_4+bC_1C_7=-C_1^{2}+C_1C_8=0$, and $C_1C_4=0$. Hence $\Gamma_{\varphi}$ is connected, a contradiction to Remark \ref{keypoint}. Therefore, $X$ admits no totally nondegenerate finite morphism. Especially, the nonsingular toric Fano $4$-folds of types $\mathcal{M}_2$ and $\mathcal{M}_4$ admit no totally nondegenerate finite morphism.

If $a=b=c=0$, then $C_3C_4=0$ as above. On the other hand, by $(5)$, we have $C_1C_3=C_1^{2}-C_1C_8=0$, and $C_1C_3=0$. Hence $\Gamma_{\varphi}$ is connected, a contradiction to Remark \ref{keypoint}. Therefore, $X$ admits no totally nondegenerate finite morphism. Especially, the nonsingular toric Fano $4$-fold of type $\mathcal{M}_1$ admits no totally nondegenerate finite morphism.

If $a=c=1$, then $C_1C_4=0$ as above. On the other hand, by $(6)$, we have $C_1C_5+C_1C_6+bC_1C_7=-C_1^{2}+C_1C_8=0$, and $C_1C_6=0$. Hence $\Gamma_{\varphi}$ is connected, a contradiction to Remark \ref{keypoint}. Therefore, $X$ admits no totally nondegenerate finite morphism. Especially, the nonsingular toric Fano $4$-fold of type $\mathcal{M}_3$ admits no totally nondegenerate finite morphism.

\bigskip

(f) ``type $\mathcal{M}_5$''\quad Let $X$ be the nonsingular projective toric $4$-fold corresponding to the fan $\Sigma$ whose primitive relations are
$$x_1+x_8=x_5,\ x_4+x_5=x_7,\ x_6+x_7=x_1,\ x_1+x_2+x_3=x_6,\ x_2+x_3+x_5=x_6+x_8$$
$$x_2+x_3+x_7=0\ \mbox{and}\ x_4+x_6+x_8=0,$$
where $\G(\Sigma)=\{x_1,\ldots,x_8\}$. $X$ is the nonsingular toric Fano $4$-fold of type $\mathcal{M}_5$. $D_1D_8=D_4D_5=D_6D_7=0$ on $X$, and by a basis $\{x_1,x_2,x_4,x_6\}$ of $N$, we have
$$(1)\ D_1-D_3+D_5+D_7=0,\ (2)\ D_2-D_3=0,$$
$$(3)\ D_4-D_5-D_8=0\ \mbox{and}\ (4)\ D_3-D_5+D_6-D_7-D_8=0$$
in ${\rm Pic}(X)$, respectively. Suppose that there exists a totally nondegenerate finite morphism $\varphi:A\rightarrow X$. By $(1)$, $(3)$ and $(4)$, we have $C_5C_8=-C_4C_5+C_5^{2}=0$ and $C_6C_8=-C_1C_8+C_8^{2}=0$, respectively. On the other hand, by Lemma \ref{conntocomp}, $C_1C_5=C_1C_7=0$. So by $(1)$, we have $C_1C_3=C_1^{2}+C_1C_5+C_1C_7=0$. Hence $\Gamma_{\varphi}$ is connected, a contradiction to Remark \ref{keypoint}. Therefore, $X$ admits no totally nondegenerate finite morphism.

\begin{Rem}\label{vksymmetric}
{\rm
Kajiwara \cite{kajiwara2}, \cite{kajiwara1} showed that the pseudo del Pezzo $4$-fold $\widetilde{V}^{4}$ admits no totally nondegenerate finite morphism similarly as above.
}

\end{Rem}

\section{The main results}\label{mainresult}

\hspace{.5cm} By Examples \ref{g1noexist} and \ref{l5noexist}, Propositions \ref{ainfdivcont} and \ref{faoverp2}, the results (a), (b), (c), (d), (e) and (f), and Remark \ref{vksymmetric}, the nonsingular toric Fano $4$-folds of types $\mathcal{B}_1$, $\mathcal{B}_2$, $\mathcal{B}_3$, $\mathcal{D}_1$, $\mathcal{D}_2$, $\mathcal{D}_3$, $\mathcal{D}_5$, $\mathcal{D}_6$, $\mathcal{D}_8$, $\mathcal{D}_9$, $\mathcal{D}_{12}$, $\mathcal{D}_{16}$, $\mathcal{G}_1$, $\mathcal{I}_4$, $\mathcal{I}_6$, $\mathcal{I}_{12}$, $\mathcal{I}_{15}$, $\mathcal{J}_2$, $\mathcal{L}_1$, $\mathcal{L}_2$, $\mathcal{L}_5$, $\mathcal{L}_{10}$, $\mathcal{L}_{12}$, $\mathcal{M}_1$, $\mathcal{M}_2$, $\mathcal{M}_3$, $\mathcal{M}_4$, $\mathcal{M}_5$ and $\widetilde{V}^{4}$ admit no totally nondegenerate embedding. Moreover, by Corollary \ref{dia2}, any nonsingular projective toric $4$-fold admits no totally nondegenerate embedding, if it is obtained by finite successions of $2$-blow-ups from one of them.

\begin{Rem}
{\rm
Kajiwara \cite{kajiwara1} and Sankaran \cite{sankaran1} showed that the nonsingular toric Fano $4$-folds of types $\mathcal{B}_5$, $\mathcal{D}_{19}$, $\mathcal{G}_{2}$ and $\mathcal{G}_6$ admit no totally nondegenerate embedding using more complicated methods (see Kajiwara \cite{kajiwara1} for types $\mathcal{D}_{19}$, $\mathcal{G}_{2}$ and $\mathcal{G}_6$, Sankaran \cite{sankaran1} for type $\mathcal{B}_5$). Since their method differs from ours, we cannot determine whether nonsingular projective toric $4$-folds obtained by finite successions of $2$-blow-ups from one of these types admit a totally nondegenerate embedding or not.
}

\end{Rem}

To describe the main result, we need the following proposition.

\begin{Prop}\label{dptimesdp}
If $X$ is a nonsingular toric Fano $4$-fold such that $X\cong X_{1}\times X_{2}$, where $X_{1}$ and $X_{2}$ are nonsingular toric del Pezzo surfaces, then there exists a totally nondegenerate embedding.

\end{Prop}

\proof
A smooth element $E_1$ in $|-K_{X_1}|$ (resp. $E_2$ in $|-K_{X_2}|$) is an elliptic curve. By an easy calculation of intersection numbers, $E_1 \times E_2\hookrightarrow X$ is obviously a totally nondegenerate embedding.\hfill q.e.d.

\begin{Rem}
{\rm
In Proposition \ref{dptimesdp}, if there exists an abelian surface embedding $A\hookrightarrow X$, then $A$ is isomorphic to the direct product of two elliptic curves as stated in the proof of Proposition \ref{dptimesdp}, by the results of Kajiwara \cite{kajiwara2} and \cite{kajiwara1}.
}

\end{Rem}

By these results and Table $1$ in Sato \cite{sato1}, we get the following:

\begin{Thm}\label{mt}
Let $X$ be a nonsingular toric Fano $4$-fold. Then, one of the following holds.

\begin{enumerate}

\item $X$ admits no totally nondegenerate embedding.

\item $X\cong {\bf P}^{4}$ or $X\cong {\bf P}^{1}\times{\bf P}^{3}$. There exists a totally nondegenerate embedding in this case $($see {\rm Horrocks-Mumford \cite{horrocks1}} and {\rm Lange \cite{lange1}}$)$.

\item $X\cong X_{1}\times X_{2}$, where $X_{1}$ and $X_{2}$ are nonsingular toric del Pezzo surfaces. There exists a totally nondegenerate embedding in this case $($see Proposition $\ref{dptimesdp})$.
\item $X$ is of one of the types $\mathcal{C}_1$, $\mathcal{C}_2$, $\mathcal{C}_3$, $\mathcal{D}_7$, $\mathcal{D}_{10}$, $\mathcal{D}_{11}$, $\mathcal{D}_{14}$, $\mathcal{D}_{17}$, $\mathcal{D}_{18}$, $\mathcal{G}_3$, $\mathcal{G}_4$, $\mathcal{G}_5$, $\mathcal{L}_{11}$, $\mathcal{L}_{13}$, $\mathcal{I}_9$, $\mathcal{Q}_{16}$, $\mathcal{U}_8$, $V^4$, $\mathcal{Z}_1$, $\mathcal{Z}_2$ and $\mathcal{W}$.

\end{enumerate}

\end{Thm}

\begin{Rem}
{\rm
For the nonsingular toric Fano $4$-fold $X$ of type $\mathcal{C}_1$, Sankaran \cite{sankaran1} showed that there exists a totally nondegenerate embedding $A\hookrightarrow X$. However, his paper seems to contain gaps unfortunately. So we do not yet know whether $X$ admits a totally nondegenerate embedding or not.
}

\end{Rem}

\section{Table of nonsingular toric Fano $4$-folds}\label{tablefano}

\hspace{.5cm} In this section, we give the table of nonsingular toric Fano $4$-folds classified in Batyrev \cite{batyrev4} and Sato \cite{sato1} with $2$-blow-up relations among them. We describe the results about totally nondegenerate embeddings obtained in the previous sections. In the third column, we show whether the nonsingular toric Fano $4$-fold admits a totally nondegenerate embedding or not. The symbol ``$\exists$'' means that there exists a totally nondegenerate embedding, while the symbol ``$\times$'' means that there does not exist a totally nondegenerate embedding. We omit a reference in the case where the noningular toric Fano $4$-fold is obtained by finite successions of $2$-blow-ups from one of the nonsingular toric Fano $4$-folds of types $\mathcal{B}_1$, $\mathcal{B}_2$, $\mathcal{B}_3$, $\mathcal{D}_1$, $\mathcal{D}_2$, $\mathcal{D}_3$, $\mathcal{D}_5$, $\mathcal{D}_6$, $\mathcal{D}_8$, $\mathcal{D}_9$, $\mathcal{D}_{12}$, $\mathcal{D}_{16}$, $\mathcal{G}_1$, $\mathcal{I}_4$, $\mathcal{I}_6$, $\mathcal{I}_{12}$, $\mathcal{I}_{15}$, $\mathcal{J}_2$, $\mathcal{L}_1$, $\mathcal{L}_2$, $\mathcal{L}_5$, $\mathcal{L}_{10}$, $\mathcal{L}_{12}$, $\mathcal{M}_1$, $\mathcal{M}_2$, $\mathcal{M}_3$, $\mathcal{M}_4$, $\mathcal{M}_5$ and $\widetilde{V}^{4}$ (see Corollary \ref{dia2}).
\setlongtables
\bigskip

\bigskip

\begin{center}
Table 1: nonsingular toric Fano $4$-folds
\end{center}

\begin{longtable}{|c||l|l||c|}

\hline
\endhead
\hline
\endfoot

\hspace{1cm} &  $2$-blow-up of\hspace{2.0cm} & embedding &notation\hspace{0.5cm} \\ \hline
(1) & none & $\exists$ (See Horrocks-Mumford \cite{horrocks1}) & ${\bf P}^{4}$ \\ \hline
(2) & none & $\times$ (See Proposition \ref{ainfdivcont}) & $\mathcal{B}_1$ \\ \hline
(3) & none & $\times$ (See Proposition \ref{ainfdivcont}) & $\mathcal{B}_2$ \\ \hline
(4) & none & $\times$ (See Proposition \ref{ainfdivcont}) & $\mathcal{B}_3$ \\ \hline
(5) & none & $\exists$ (See Lange \cite{lange1}) & $\mathcal{B}_4$ \\ \hline
(6) & ${\bf P}^{4}$ & $\times$ (See Sankaran \cite{sankaran1}) & $\mathcal{B}_5$ \\ \hline
(7) & none & unknown & $\mathcal{C}_1$ \\ \hline
(8) & none & unknown & $\mathcal{C}_2$ \\ \hline
(9) & none & unknown & $\mathcal{C}_3$ \\ \hline
(10) & none & $\exists$ (See Proposition \ref{dptimesdp}) & $\mathcal{C}_4$ \\ \hline
(11) & $\mathcal{B}_1,\ \mathcal{B}_2$ & $\times$ & $\mathcal{E}_1$ \\ \hline
(12) & $\mathcal{B}_2,\ \mathcal{B}_3$ & $\times$ & $\mathcal{E}_2$ \\ \hline
(13) & $\mathcal{B}_3,\ \mathcal{B}_4$ & $\times$ & $\mathcal{E}_3$ \\ \hline
(14) & none & $\times$ (See Proposition \ref{faoverp2}) & $\mathcal{D}_1$ \\ \hline
(15) & $\mathcal{C}_1$ & $\times$ (See Proposition \ref{faoverp2}) & $\mathcal{D}_2$ \\ \hline
(16) & none & $\times$ (See Proposition \ref{faoverp2}) & $\mathcal{D}_3$ \\ \hline
(17) & $\mathcal{B}_2$ & $\times$ & $\mathcal{D}_4$ \\ \hline
(18) & none & $\times$ (See Proposition \ref{faoverp2}) & $\mathcal{D}_5$ \\ \hline
(19) & $\mathcal{C}_3$ & $\times$ (See Proposition \ref{faoverp2}) & $\mathcal{D}_6$ \\ \hline
(20) & none & unknown & $\mathcal{D}_7$ \\ \hline
(21) & $\mathcal{C}_2$ & $\times$ (See Proposition \ref{faoverp2}) & $\mathcal{D}_8$ \\ \hline
(22) & none & $\times$ (See Proposition \ref{faoverp2}) & $\mathcal{D}_9$ \\ \hline
(23) & $\mathcal{B}_5$ & unknown & $\mathcal{D}_{10}$ \\ \hline
(24) & $\mathcal{B}_5,\ \mathcal{C}_2$ & unknown & $\mathcal{D}_{11}$ \\ \hline
(25) & none & $\times$ (See Proposition \ref{faoverp2}) & $\mathcal{D}_{12}$ \\ \hline
(26) & none & $\exists$ (See Proposition \ref{dptimesdp}) & $\mathcal{D}_{13}$ \\ \hline
(27) & $\mathcal{B}_4$ & unknown & $\mathcal{D}_{14}$ \\ \hline
(28) & $\mathcal{C}_4$ & $\exists$ (See Proposition \ref{dptimesdp}) & $\mathcal{D}_{15}$ \\ \hline
(29) & $\mathcal{C}_3$ & $\times$ (See Proposition \ref{faoverp2}) & $\mathcal{D}_{16}$ \\ \hline
(30) & $\mathcal{B}_5$ & unknown & $\mathcal{D}_{17}$ \\ \hline
(31) & $\mathcal{C}_1$ & unknown & $\mathcal{D}_{18}$ \\ \hline
(32) & $\mathcal{C}_2$ & $\times$ (See Kajiwara \cite{kajiwara1}) & $\mathcal{D}_{19}$ \\ \hline
(33) & none & $\times$ (See Example \ref{g1noexist}) & $\mathcal{G}_1$ \\ \hline
(34) & $\mathcal{C}_2$ & $\times$ (See Kajiwara \cite{kajiwara1}) & $\mathcal{G}_2$ \\ \hline
(35) & none & unknown & $\mathcal{G}_3$ \\ \hline
(36) & $\mathcal{C}_2$ & unknown & $\mathcal{G}_4$ \\ \hline
(37) & $\mathcal{C}_3$ & unknown & $\mathcal{G}_5$ \\ \hline
(38) & $\mathcal{C}_4$ & $\times$ (See Kajiwara \cite{kajiwara1}) & $\mathcal{G}_6$ \\ \hline
(39) & $\mathcal{D}_2$ & $\times$ & $\mathcal{H}_1$ \\ \hline
(40) & $\mathcal{D}_3$ & $\times$ & $\mathcal{H}_2$ \\ \hline
(41) & $\mathcal{D}_1,\ \mathcal{D}_5$ & $\times$ & $\mathcal{H}_3$ \\ \hline
(42) & $\mathcal{D}_8,\ \mathcal{D}_9$ & $\times$ & $\mathcal{H}_4$ \\ \hline
(43) & $\mathcal{D}_6,\ \mathcal{D}_{12},\ \mathcal{D}_{16}$ & $\times$ & $\mathcal{H}_5$ \\ \hline
(44) & $\mathcal{D}_3,\ \mathcal{D}_9$ & $\times$ & $\mathcal{H}_6$ \\ \hline
(45) & $\mathcal{D}_2,\ \mathcal{D}_5,\ \mathcal{D}_{18}$ & $\times$ & $\mathcal{H}_7$ \\ \hline
(46) & $\mathcal{D}_{13},\ \mathcal{D}_{15}$ & $\exists$ (See Proposition \ref{dptimesdp}) & $\mathcal{H}_8$ \\ \hline
(47) & $\mathcal{D}_8,\ \mathcal{D}_{12},\ \mathcal{D}_{19}$ & $\times$ & $\mathcal{H}_9$ \\ \hline
(48) & $\mathcal{D}_9,\ \mathcal{D}_{16}$ & $\times$ & $\mathcal{H}_{10}$ \\ \hline
(49) & none & $\times$ (See $(c)$ in Section \ref{secex}) & $\mathcal{L}_1$ \\ \hline
(50) & $\mathcal{D}_7$ & $\times$ (See $(c)$ in Section \ref{secex}) & $\mathcal{L}_2$ \\ \hline
(51) & $\mathcal{D}_6$ & $\times$ & $\mathcal{L}_3$ \\ \hline
(52) & $\mathcal{D}_8,\ \mathcal{D}_{10},\ \mathcal{D}_{11}$ & $\times$ & $\mathcal{L}_4$ \\ \hline
(53) & none & $\times$ (See Example \ref{l5noexist}) & $\mathcal{L}_5$ \\ \hline
(54) & $\mathcal{D}_{12},\ \mathcal{D}_{14}$ & $\times$ & $\mathcal{L}_6$ \\ \hline
(55) & $\mathcal{D}_{15}$ & $\exists$ (See Proposition \ref{dptimesdp}) & $\mathcal{L}_7$ \\ \hline
(56) & none & $\exists$ (See Proposition \ref{dptimesdp}) & $\mathcal{L}_8$ \\ \hline
(57) & $\mathcal{D}_{13}$ & $\exists$ (See Proposition \ref{dptimesdp}) & $\mathcal{L}_9$ \\ \hline
(58) & $\mathcal{D}_{10},\ \mathcal{D}_{17}$ & $\times$ (See $(c)$ in Section \ref{secex}) & $\mathcal{L}_{10}$ \\ \hline
(59) & $\mathcal{D}_{14}$ & unknown & $\mathcal{L}_{11}$ \\ \hline
(60) & $\mathcal{D}_{11},\ \mathcal{D}_{17},\ \mathcal{D}_{19}$ & $\times$ (See $(d)$ in Section \ref{secex}) & $\mathcal{L}_{12}$ \\ \hline
(61) & $\mathcal{D}_7$ & unknown & $\mathcal{L}_{13}$ \\ \hline
(62) & $\mathcal{D}_4$ & $\times$ & $\mathcal{I}_1$ \\ \hline
(63) & $\mathcal{D}_1,\ \mathcal{D}_6$ & $\times$ & $\mathcal{I}_2$ \\ \hline
(64) & $\mathcal{D}_3,\ \mathcal{D}_8$ & $\times$ & $\mathcal{I}_3$ \\ \hline
(65) & $\mathcal{D}_{10}$ & $\times$ (See $(a)$ in Section \ref{secex}) & $\mathcal{I}_4$ \\ \hline
(66) & $\mathcal{E}_2,\ \mathcal{D}_4,\ \mathcal{D}_{10}$ & $\times$ & $\mathcal{I}_5$ \\ \hline
(67) & $\mathcal{D}_{10}$ & $\times$ (See $(a)$ in Section \ref{secex}) & $\mathcal{I}_6$ \\ \hline
(68) & $\mathcal{D}_5,\ \mathcal{D}_{12}$ & $\times$ & $\mathcal{I}_7$ \\ \hline
(69) & $\mathcal{D}_8,\ \mathcal{D}_{16}$,\ $\mathcal{G}_4$ & $\times$ & $\mathcal{I}_8$ \\ \hline
(70) & $\mathcal{D}_{14}$ & unknown & $\mathcal{I}_9$ \\ \hline
(71) & $\mathcal{D}_6,\ \mathcal{D}_{15},\ \mathcal{G}_5$ & $\times$ & $\mathcal{I}_{10}$ \\ \hline
(72) & $\mathcal{D}_9,\ \mathcal{D}_{12}$ & $\times$ & $\mathcal{I}_{11}$ \\ \hline
(73) & $\mathcal{D}_{15},\ \mathcal{D}_{19},\ \mathcal{G}_6$ & $\times$ (See $(a)$ in Section \ref{secex}) & $\mathcal{I}_{12}$ \\ \hline
(74) & $\mathcal{D}_{12},\ \mathcal{D}_{13}$ & $\times$ & $\mathcal{I}_{13}$ \\ \hline
(75) & $\mathcal{E}_3,\ \mathcal{D}_{10},\ \mathcal{D}_{14}$ & $\times$ & $\mathcal{I}_{14}$ \\ \hline
(76) & $\mathcal{D}_{18},\ \mathcal{D}_{19},\ \mathcal{G}_2$ & $\times$ (See $(a)$ in Section \ref{secex}) & $\mathcal{I}_{15}$ \\ \hline
(77) & none & $\times$ (See $(e)$ in Section \ref{secex}) & $\mathcal{M}_1$ \\ \hline
(78) & none & $\times$ (See $(e)$ in Section \ref{secex}) & $\mathcal{M}_2$ \\ \hline
(79) & $\mathcal{G}_3,\ \mathcal{G}_5$ & $\times$ (See $(e)$ in Section \ref{secex}) & $\mathcal{M}_3$ \\ \hline
(80) & $\mathcal{G}_3$ & $\times$ (See $(e)$ in Section \ref{secex}) & $\mathcal{M}_4$ \\ \hline
(81) & $\mathcal{G}_4,\ \mathcal{G}_6$ & $\times$ (See $(f)$ in Section \ref{secex}) & $\mathcal{M}_5$ \\ \hline
(82) & $\mathcal{G}_1,\ \mathcal{G}_3$ & $\times$ & $\mathcal{J}_1$ \\ \hline
(83) & $\mathcal{G}_3$ & $\times$ (See $(b)$ in Section \ref{secex}) & $\mathcal{J}_2$ \\ \hline
(84) & $\mathcal{L}_2$ & $\times$ & $\mathcal{Q}_1$ \\ \hline
(85) & $\mathcal{H}_4,\ \mathcal{L}_4$ & $\times$ & $\mathcal{Q}_2$ \\ \hline
(86) & $\mathcal{L}_1,\ \mathcal{L}_5$ & $\times$ & $\mathcal{Q}_3$ \\ \hline
(87) & $\mathcal{L}_3$ & $\times$ & $\mathcal{Q}_4$ \\ \hline
(88) & $\mathcal{H}_5,\ \mathcal{L}_3,\ \mathcal{L}_6$ & $\times$ & $\mathcal{Q}_5$ \\ \hline
(89) & $\mathcal{L}_6$ & $\times$ & $\mathcal{Q}_6$ \\ \hline
(90) & $\mathcal{L}_7$ & $\times$ & $\mathcal{Q}_7$ \\ \hline
(91) & $\mathcal{L}_5,\ \mathcal{L}_9$ & $\times$ & $\mathcal{Q}_8$ \\ \hline
(92) & $\mathcal{L}_3,\ \mathcal{L}_7,\ \mathcal{I}_{10}$ & $\times$ & $\mathcal{Q}_9$ \\ \hline
(93) & $\mathcal{H}_8,\ \mathcal{L}_7,\ \mathcal{L}_9$ & $\exists$ (See Proposition \ref{dptimesdp}) & $\mathcal{Q}_{10}$ \\ \hline
(94) & $\mathcal{L}_8,\ \mathcal{L}_9$ & $\exists$ (See Proposition \ref{dptimesdp}) & $\mathcal{Q}_{11}$ \\ \hline
(95) & $\mathcal{L}_{10},\ \mathcal{L}_{12},\ \mathcal{I}_6$ & $\times$ & $\mathcal{Q}_{12}$ \\ \hline
(96) & $\mathcal{L}_2,\ \mathcal{L}_5,\ \mathcal{L}_{13}$ & $\times$ & $\mathcal{Q}_{13}$ \\ \hline
(97) & $\mathcal{H}_9,\ \mathcal{L}_4,\ \mathcal{L}_6,\ \mathcal{L}_{12},\ \mathcal{I}_{14}$ & $\times$ & $\mathcal{Q}_{14}$ \\ \hline
(98) & $\mathcal{L}_6,\ \mathcal{L}_9,\ \mathcal{L}_{11},\ \mathcal{I}_{13}$ & $\times$ & $\mathcal{Q}_{15}$ \\ \hline
(99) & $\mathcal{L}_{11},\ \mathcal{L}_{13},\ \mathcal{I}_9$ & unknown & $\mathcal{Q}_{16}$ \\ \hline
(100) & $\mathcal{L}_7,\ \mathcal{L}_{12},\ \mathcal{I}_{12}$ & $\times$ & $\mathcal{Q}_{17}$ \\ \hline
(101) & $\mathcal{H}_1,\ \mathcal{H}_3,\ \mathcal{H}_7$ & $\times$ & $\mathcal{K}_1$ \\ \hline
(102) & $\mathcal{H}_2,\ \mathcal{H}_6,\ \mathcal{H}_{10}$ & $\times$ & $\mathcal{K}_2$ \\ \hline
(103) & $\mathcal{H}_4,\ \mathcal{H}_5,\ \mathcal{H}_9$ & $\times$ & $\mathcal{K}_3$ \\ \hline
(104) & $\mathcal{H}_8$ & $\exists$ (See Proposition \ref{dptimesdp}) & $\mathcal{K}_4$ \\ \hline
(105) & $\mathcal{M}_3$ & $\times$ & $\mathcal{R}_1$ \\ \hline
(106) & $\mathcal{M}_2,\ \mathcal{M}_4$ & $\times$ & $\mathcal{R}_2$ \\ \hline
(107) & $\mathcal{M}_1,\ \mathcal{M}_4$ & $\times$ & $\mathcal{R}_3$ \\ \hline
(108) & $\mathcal{I}_{11},\ \mathcal{I}_{13}$ & $\times$ & $\mathcal{P}$ \\ \hline
(109) & $\mathcal{Q}_1,\ \mathcal{Q}_3,\ \mathcal{Q}_{13}$ & $\times$ & $\mathcal{U}_1$ \\ \hline
(110) & $\mathcal{Q}_2,\ \mathcal{Q}_5,\ \mathcal{Q}_{14},\ \mathcal{K}_3$ & $\times$ & $\mathcal{U}_2$ \\ \hline
(111) & $\mathcal{Q}_4,\ \mathcal{Q}_9$ & $\times$ & $\mathcal{U}_3$ \\ \hline
(112) & $\mathcal{Q}_{10},\ \mathcal{K}_4$ & $\times$ & $\mathcal{U}_4$ \\ \hline
(113) & $\mathcal{Q}_{11}$ & $\exists$ (See Proposition \ref{dptimesdp}) & $\mathcal{U}_5$ \\ \hline
(114) & $\mathcal{Q}_6,\ \mathcal{Q}_8,\ \mathcal{Q}_{15}$ & $\times$ & $\mathcal{U}_6$ \\ \hline
(115) & $\mathcal{Q}_7,\ \mathcal{Q}_{12},\ \mathcal{Q}_{17}$ & $\times$ & $\mathcal{U}_7$ \\ \hline
(116) & $\mathcal{Q}_{16}$ & unknown & $\mathcal{U}_8$ \\ \hline
(117) & none & $\times$ (See Kajiwara \cite{kajiwara2} and \cite{kajiwara1}) & ${\widetilde{V}}^{4}$ \\ \hline
(118) & none & unknown & $V^{4}$ \\ \hline
(119) & $\mathcal{Q}_{10},\ \mathcal{Q}_{11}$ & $\exists$ (See Proposition \ref{dptimesdp}) & $S_2\times S_2$ \\ \hline
(120) & $\mathcal{U}_4,\ \mathcal{U}_5,\ S_2\times S_2$ & $\exists$ (See Proposition \ref{dptimesdp}) & $S_2\times S_3$ \\ \hline
(121) & $S_2\times S_3$ & $\exists$ (See Proposition \ref{dptimesdp}) & $S_3\times S_3$ \\ \hline
(122) & $\mathcal{G}_6$ & unknown & $\mathcal{Z}_1$ \\ \hline
(123) & $\mathcal{G}_4$ & unknown & $\mathcal{Z}_2$ \\ \hline
(124) & $\mathcal{Z}_1$ & unknown & $\mathcal{W}$ \\ \hline

\end{longtable}

\bigskip

\begin{flushleft}
\begin{sc}
Department of Mathematics,\\
Tokyo Institute of Technology,\\
2-12-1 Oh-okayama, Meguro-ku, Tokyo 152-8551, Japan.
\end{sc}

\medskip
{\it E-mail address}: hirosato@math.titech.ac.jp
\end{flushleft}

\end{document}